\documentclass[12pt]{article}
\usepackage{amssymb}
\usepackage{amsthm}
\usepackage{lineno}
\usepackage{amsfonts}
\usepackage{amsmath}

\newtheorem{thm}{Theorem}

\newtheorem{lemma}{Lemma}

\begin{document}
\title{ Généralisation du Théorème de Zeckendorf}
\author{\textbf{Rachid CHERGUI} \thanks{%
The research is partially supported by ATN laboratory of USTHB University.%
} \\
{\small USTHB, Faculty of mathematics }\\
[-0.8ex] {\small P.B. 32, El Alia, 16111, Bab Ezzouar, Algeria}\\
[-0.8ex] {\small \texttt{rchergui@usthb.dz 
}}
}
\date{ {\small Mathematics Subject Classifications: 05A10,
05A17, 11B39, 11B50.}}
\maketitle
\begin{abstract}
  On  consiste a présenter d’abord le Théorème de Zeckendorf
avec ces deux versions Fibonacci et Luca .\\ 
Dans ce papier on obtient un résultas   sur la généralisation du théorème Zeckendorf  pour les  nombres de Fibonacci généralisé  (multibonacci). De tels résultats trouvent des application dans la théorie du codage.
\end{abstract}
\maketitle
\section{Introduction}
Depuis l’existence de la suite de Fibonacci et lucas, les chercheurs ne s’arrêtent pas a 
trouver des applications de cette suite dans des différents domaines.
Le théorème de Zeckendorf, dénommé ainsi d'après le mathématicien belge Édouard Zeckendorf, est un théorème de théorie additive des nombres qui garantit que tout entier naturel N peut être représenté, de manière unique, comme somme de nombres de Fibonacci distincts et non consécutifs. Cette représentation est appelée la représentation de Zeckendorf de $ N$.\\
Rachid chergui a étudier  sur l'arithmétique de Zeckendorf. Dans ce papier on s'intéresse a sa généralisation.
\newpage
\section{ Théorème de Zeckendorf pour les nombres de Fibonacci}
Le Théorème de Zeckendorf, nommé de son fondateur belge Edouard Zeckendorf, est un théorème pour la représentation des entiers comme somme de nombres de Fibonacci.
\begin{thm}
 Chaque entier positif peut être représenté d'une manière unique comme somme de nombres de Fibonacci de tel sorte qu'il n'existe pas deux nombres de Fibonacci consécutifs. Autrement dit, si $n$ est un entier positif

$$
n=\sum_{r=1}^{\infty} e_{r} F_{r}
$$
où $e_{r} \in\{0,1\}, e_{r}=1 \Rightarrow e_{r+1}=0$ et pour un nombre fini de coefficients $e_{r}$ sont non nuls.
\end{thm}
Preuve. le Théorème de Zeckendorf a deux parties :

Existence : chaque entier positif $n$ a une représentation de Zeckendorf.

Unicité : chaque entier positif n'a pas deux représentations différentes de Zeckendorf.\\
La première partie du Théorème de Zeckendorf(existence) peut être prouvée par induction. Pour $n=1,2,3$ elle est vraie, pour $n=4$ on a $4=3+1$. Maintenant, supposons que chaque entier a une représentation de Zeckendorf. Si $k+1$ est un nombre de Fibonacci alors c'est fini, sinon il existe $j$ tel que $F_{j}<k+1<F_{j+1}$. Considérons $a=k+1-F_{j}$.

Comme $a<k, a$ a une représentation de Zeckendorf; de plus $F_{j}+a<F_{j+1}$ donc $a<F_{j-1}$ donc la représentation de Zeckendorf de $a$ ne contient pas $F_{j-1}$. Alors $k+1$ peut être représenté comme $a+F_{j}$. De plus, il est clair que chaque représentation de Zeckendorf correspond à un seul entier. La seconde partie du Théorème de Zeckendorf (unicité) requiert le lemme suivant :
\begin{lemma}
 La somme des éléments d'un ensemble non vide de nombres de Fibonacci distincts, non-consécutifs de plus grand nombre $F_{j}$ est strictement inférieur à $F_{j+1}$.
\end{lemma}
Preuve. Maintenant, on prend deux ensembles de nombres de Fibonacci distincts, non-consécutifs $S$ et $T$ qui ont la même somme. Éliminons les nombres communs pour former des ensembles $S_{\prime}$ et $T \prime$ qui n'ont aucun nombre commun. On veut montrer que $S \prime$ et $T \prime$ sont vides i.e $S=T$.

Premièrement, on montre qu'au moins un des ensembles $S^{\prime}$ et $T \prime$ est vide. Supposons le contraire, soit $F_{s}$ le plus grand nombre de $S^{\prime}$ et $F_{t}$ le plus grand nombre de $T \prime$ sans perte de
généralité, supposons que $F_{s}<F_{t}$. Alors par Lemme 13, la somme de $S^{\prime}$ est strictement inférieur à $F_{s+1}$, et donc strictement inférieur à $F_{t}$, mais il est claire que la somme de $T \prime$ est au moins $F_{t}$.

Cela veut dire que $S$ ' et $T \prime$ ne peuvent pas avoir la même somme, et alors $S$ et $T$ ne peuvent pas avoir la même somme.

Donc notre supposition est fausse.

Si $S^{\prime}$ est vide et $T^{\prime}$ est non vide alors $S$ est un sous ensemble de $T$, et alors $S$ et $T$ ne peuvent pas avoir la même somme. Similairement on peut éliminer le cas ou $S$ ' est non vide et $T \prime$ est vide. Le seul cas qui reste est que $S \prime$ et $T \prime$ sont vides, donc $S=T$.

On conclut que quelque soient deux représentations de Zeckendorf qui ont la même somme doivent être identiques.
\section{ Théorème de Zeckendorf pour les nombres de Lucas}
Il est connu que les nombres de Lucas sont complets\cite{1} au sens que chaque entier positif peut être représenté comme somme de nombres distincts de Lucas. En général, ses représentation ne sont pas uniques par exemple $4=L_{3}=L_{1}+L_{2}, 12=L_{1}+L_{3}+L_{4}=L_{0}+L_{2}+L_{4}$.

Avant d'énoncer le Théorème, certains Lemmes sont utiles :
\begin{lemma}
 $ L_{n}-1=L_{n-1}+L_{n-3}+L_{n-5}+\ldots+L_{1,2}(n)$, pour $n \geq 2$ où

$$
L_{1,2}(n)=\left\{\begin{array}{l}
2 L_{1} \text { si } n \text { pair } \\
L_{2} \text { si } n \text { impair }
\end{array}\right.
$$
\end{lemma}
Preuve. on va monter ce Lemme par récurrence sur $n$ :
Si $n$ est pair, posons $n=2 k$
on a $L_{2}-1=3-1=2=2 L_{1}$ vérifie le lemme.
Supposons que l'égalité est vraie pour l'ordre $n$ et montrons qu'elle est vraie pour l'ordre $n+1$. On a
$$
L_{2 k}-1=L_{2 k-1}+L_{2 k-3}+L_{2 k-5}+\ldots+L_{3}+2 L_{1}
$$

on rajoutant $L_{2 k+1}$ des deux cotés on aura
$$
L_{2 k+1}+L_{2 k}-1=L_{2 k+1}+L_{2 k-1}+L_{2 k-3}+\ldots+L_{3}+2 L_{1}
$$

or $L_{2 k+1}+L_{2 k}=L_{2 k+2}$ donc
$$
L_{2 k+2}=L_{2 k+1}+L_{2 k-1}+L_{2 k-3}+\ldots+L_{3}+2 L_{1}
$$
ce ci est égale à
$$
L_{n+2}=L_{(n+2)-1}+L_{(n+2)-3}+L_{(n+2)-5}+\ldots+L_{3}+2 L_{1}
$$
et on a bien l'égalité du lemme.

si $\mathrm{n}$ est impaire, posons $n=2 k+1$.

on a $L_{3}-1=4-1=3=L_{2}$, vérifie le lemme.

Supposons que l'égalité est vraie pour l'ordre $n$ et montrons qu'elle est vraie pour l'ordre $n+1$. On a

$$
L_{2 k+1}-1=L_{2 k}+L_{2 k-2}+L_{2 k-4}+\ldots+L_{4}+L_{2}
$$

on rajoutant $L_{2 k+2}$ des deux cotés on aura

$$
L_{2 k+2}+L_{2 k+1}-1=L_{2 k+2}+L_{2 k}+L_{2 k-2}+L_{2 k-4}+\ldots+L_{4}+L_{2}
$$

or $L_{2 k+2}+L_{2 k+1}=L_{2 k+3}$, donc

$$
L_{2 k+3}=L_{2 k+2}+L_{2 k}+L_{2 k-2}+L_{2 k-4}+\ldots+L_{4}+L_{2}
$$

en remplaçons $n$ par $2 k+1$ dans l'égalité précédente on aura le résultat

$$
L_{n+2}=L_{(n+2)-1}+L_{(n+2)-3}+L_{(n+2)-5}+\ldots+L_{4}+L_{2}
$$

Donc le lemme est vrai dans les deux cas.

\begin{lemma}
$ L_{n+2}=1+\sum_{i=0}^{n} L_{i}$ pour $n \geq 0$.
\end{lemma}
Preuve. On va démontrer ce Lemme par récurrence sur $n$.

pour $n=0$ on a $L_{2}=1+L_{1}=1+2=3$. L'égalité est vraie.

Supposons qu'elle est vraie pour $n$ est montrons qu'elle est vraie pour $n+1$.

on rajoutant $L_{n+1}$ des deux cotés de l'égalité on aura

$$
L_{n+2}+L_{n+1}=1+\sum_{i=0}^{n} L_{i}+L_{n+1}
$$

ce ci est égale à

$$
L_{n+3}=L_{(n+1)+2}=1+\sum_{i=0}^{n+1} L_{i}
$$

ce qui prouve le Lemme.
\begin{thm}
 Soit $n$ un entier positif satisfaisant $0 \leq n \leq L_{k}$ pour $k \geq 1$, alors
$$
n=\sum_{0}^{k-1} \alpha_{i} L_{i}
$$
où $\alpha_{i} \in\{0,1\}$ tels que
$$
\left\{\begin{array}{l}
\alpha_{i} \alpha_{i+1}=0 \text { pour } i \geq 0 \\
\alpha_{0} \alpha_{2}=0
\end{array}\right.
$$
Cette représentation est unique.
\end{thm}
Preuve. La preuve va ce décomposer en deux parties : existence et unicité de la représentation.

Existence : Supposons que $n$ a aussi cette représentation $n=\sum_{0}^{k-1} \gamma_{i} L_{i}$ avec $\gamma_{i} \in\{0,1\}$, $\gamma_{i} \gamma_{i+1}=0$ et $\gamma_{0} \gamma_{2}=0$.

Supposons pour une preuve par contradiction que les deux représentations ne sont pas identiques,

$$
\sum_{0}^{\infty}\left|\gamma_{i}-\alpha_{i}\right| \neq 0
$$

alors, soit $k$ la plus grande valeur de $i$ tel que $\gamma_{i} \neq \alpha_{i}$. Plus précisément $k \geq 2$, et comme $\gamma_{k} \neq \alpha_{k}$, on peut supposer, sans perte de généralité, que $\alpha_{k}=1, \gamma_{k}=0$.

Pour un $m \leq n$,

$$
m=\sum_{0}^{k} \alpha_{i} L_{i}=\sum_{0}^{k-1} \gamma_{i} L_{i}, \text { avec } \alpha_{k}=1
$$

Alors

$$
\sum_{0}^{k} \alpha_{i} L_{i} \geq L_{k}
$$

à partir des contraintes sur les coefficients $\left(\gamma_{i}\right)$,

$$
\sum_{0}^{k-1} \gamma_{i} L_{i} \leq L_{k-1}+L_{k-3}+\ldots+L_{1,2}(k)=L_{k}-1
$$

Ainsi $m \geq L_{k}$ tant que $m \leq L_{k}-1$, contradiction.

Unicité : Supposons que $n$ à deux représentations

\begin{equation}
n=\sum_{0}^{k} \beta_{i} L_{i}=\sum_{0}^{k-1} \gamma_{i} L_{i}
\end{equation}

où $\beta_{i}, \gamma_{i} \in\{0,1\}$ tel que $\beta_{k}=\gamma_{m}=1, \beta_{i}+\beta_{i+1} \neq 0$. Pour $0 \leq i \leq k-2$ :

$$
\beta_{0}+\beta_{2} \neq 0, \gamma_{i}+\gamma_{i+1} \neq 0
$$

et pour $0 \leq i \leq m-2$

$$
\gamma_{0}+\gamma_{2} \neq 0
$$

Sans perte de généralité, on prend $m \geq k \geq 2$. Si $m>k$ alors la représentation à droite en (1) avec les contraintes de constantes implique

$$
n \geq\left\{\begin{array}{l}
L_{m}+L_{m-2}+\ldots+L_{2}+L_{1}=L_{m+1} \geq L_{k+2} \quad \text { (m pair) } \\
L_{m}+L_{m-2}+\ldots+L_{3}+L_{1}+L_{0}=L_{m+1} \geq L_{k+2} \quad \text { (m impair) }
\end{array}\right.
$$

mais

$$
n=\sum_{0}^{k} \beta_{i} L_{i} \leq \sum_{0}^{k} L_{i}=L_{k+2}-1
$$

contradiction. Donc $m=k$ dans (1).
\section{ Méthode de la décomposition de Zeckendorf }
Pour décomposer un entier $x$ de la forme de Zeckendorf $x=\sum_{r=1}^{\infty} e_{r} F_{r}$ il suffit de suivre les étapes suivantes :

\begin{enumerate}
  \item Trouver le plus grand nombre de Fibonacci $F_{r} \leq x$;

  \item Effectuer la soustraction $X=x-F_{r}$; affecter un '1' à $e_{r}$ et conserver ce coefficient;

  \item Affecter $X$ à $x$ et répéter les étapes 1 et 2 jusqu'à avoir un $X$ nul;

  \item Affecter des 0 aux $e_{i}$ où $0<i<r$ et $e_{i} \neq 1$.

\end{enumerate}

Le résultat de cette décomposition est : un vecteur de $r$ éléments qui contient les coefficients $e_{r}$ de la décomposition.

\section{Résultats obtenus.\\
Généralisation du Théorème de Zeckendorf}
Pour généraliser le Théorème de Zeckendorf, ces deux Lemmes sont essentiels dans la preuve .

\begin{lemma}

$$
\sum_{r=0}^{k} 2^{r}=2^{k+1}-1
$$
\end{lemma}
\begin{lemma}
 Soit $\left(g_{n}\right)$ une suite de Fibonacci généralisé d'ordre $k$, alors pour tout entier $i, g_{i+1}^{(k)}<$ $2 g_{i}^{(k)}$.
\end{lemma}
Preuve. Par la récurrence de Fibonacci on a :
$$
2 g_{i}^{(k)}=g_{i}^{(k)}+g_{i-1}^{(k)}+g_{i-2}^{(k)}+\ldots+g_{i-(k-1)}^{(k)}+g_{i-k}^{(k)}
$$
Ceci implique que
$$
2 g_{i}^{(k)}=g_{i+1}^{(k)}+g_{i-k}^{(k)}
$$
aussi
$$
g_{i+1}^{(k)}<2 g_{i+1}^{(k)}
$$
Maintenant, On va prouver un théorème analogue au tour de la représentation de Zeckendorf pour les suites de Fibonacci d'ordre $k$.
\begin{thm}
 Considérons l'ensemble des suites recurrentes linéaires d'ordre $k\left(U_{n}^{(k)}\right)_{n}$ satisfaisant la récurrence de Fibonacci généralisée d'ordre $k$ telles que $U_{0}^{(k)}<U_{1}^{(k)}<\ldots<U_{k-2}^{(k)}<$ $U_{k-1}^{(k)}$. Dans cet ensemble il existe une unique suite recurrente linéaire $\left(U_{n}\right)$ d'ordre $k$ généralisé tel que pour tout entier positif $n$ on a une unique représentation de Zeckendorf d'ordre $k$.
\end{thm}
Preuve. En premier temps, on va démontrer par récurrence l'existence d'une unique représentation de Zeckendorf d'ordre $k$ pour tout entier positif $m$.
Il est claire que $\forall m \in \mathbb{N}, 0 \leq m \leq 2^{k}-1$.

La représentation des entiers $x$ tel que $0 \leq m \leq 2^{k}-2$ correspond à la représentation en binaire de ces nombres. Par le lemme1 $2^{0}+2^{1}+\ldots+2^{k-1}=2^{k}-1$, donc la représentation de $2^{k}-1$ est $g_{k}^{(k)}$.

Supposons que pour tout $j$ tel que $0 \leq j \leq x$ il existe une unique représentation de Zeckendorf d'ordere $k$ de $x$. Soit $i$ le plus grand entier tel que $g_{i}^{(k)} \leq x$. On sait que ce $i$
existe car les termes de $\left(g_{n}\right)$ sont croissantes. Soit $x-g_{i}^{(k)}=g_{a_{1}}^{(k)}+g_{a_{2}}^{(k)}+\ldots+g_{a_{r}}^{(k)}$ une unique représentation de $x-g_{i}^{(k)}$, alors il existe une représentation de $x=g_{i}^{(k)}+g_{a_{1}}^{(k)}+g_{a_{2}}^{(k)}+\ldots+g_{a_{r}}^{(k)}$. Cette représentation des de Zeckendorf d'ordre $k$ si $r<k-1$.
Maintenant, supposons que $r \geq k-1$. On peut poser $a_{1}>a_{2}>\ldots>a_{r}$. On va montrer que $i>a_{1}$. Si $i<a_{1}$ alors $a_{1}$ est le plus grand entier tel que $g_{a_{1}}^{(k)} \leq x$ contradiction car $i$ est le plus grand entier tel que $g_{a_{1}}^{(k)} \leq x$.
Supposons que $i=a_{1}$, alors $x=g_{i}^{(k)}+g_{a_{1}}^{(k)}+g_{a_{2}}^{(k)}+\ldots+g_{a_{r}}^{(k)}$, le lemme 2 nous dis que cette somme est plus grande que $g_{i}^{(k)}$ donc $g_{i+1}^{(k)}>x$. Alors $g_{i}^{(k)}+g_{a_{1}}^{(k)}+g_{a_{2}}^{(k)}+\ldots+g_{a_{r}}^{(k)}$ est une représentation de Zeckendorf d'ordre $k$ de $x$.
Supposons qu'il existe une autre représentation de Zeckendorf d'ordre $k$ de $x$ qui ne contient pas $g_{i}$, alors la valeur maximale de cette représentation est
$$
F(i)=\sum_{\substack{1 \leq j \leq i \\ k \nmid j}} g_{i-j}^{(k)}
$$
On utilise une récurrence sur $i$ pour démontrer $F(i)=g_{i}-1<x$ pour tout $i$. Par le lemme 1 , on sait que c'est vrai pour le cas de la basse $i \leq k-1$. Supposons que $i \geq k$ et que $F(l)=g_{l}^{(k)}-1$ pour tout $l<i$. En particulier, $F_{i-k}-1$,
$$
\begin{aligned}
F_{i} & =g_{i}^{(k)}+g_{i-1}^{(k)}+g_{i-2}^{(k)}+\ldots+g_{i-(k-1)}^{(k)}+F(i-k) \\
& =g_{i}^{(k)}-g_{i-k}^{(k)}+\left(g_{i-k}^{(k)}-1\right) \\
& =g_{i}^{(k)}-1<x .
\end{aligned}
$$
donc la représentation de Zeckendorf d'ordre $k$ est unique.
On montre maintenant l'unicité de $\left(g_{n}\right)$.
Supposons pour un ordre $k$ arbitraire, il existe une autre suite de Fibonacci généralisé $\left(h_{n}\right)$ avec $h_{0}<h_{1}<\ldots<h_{k-1}$ tel que chaque entier positif $x$ à une unique représentation de Zeckendorf d'ordre $k$ qui respecte $\left(h_{n}\right)$. Alors il reste au moins un entier $q$ dans l'intervalle $[0, k-1]$ tel que $g_{q} \neq h_{q}$. Si $h_{q}<g_{q}$ alors $h_{q}$ a plus qu'une représentation de Zeckendorf d'ordre $k$ dans $\left(h_{n}\right)$. Les deux représentations sont $h_{q}$ et la représentation en binaire de $h_{p}$ dans $h_{0}, h_{1}, \ldots, h_{q-1}$. Si $h_{p}>g_{q}$ alors $2^{q}$ n'a pas une représentation de Zeckendorf d'ordre $k$ dans $\left(h_{n}\right)$.
Comme $h_{q}>g_{q}^{(k)}=2^{q} ; 2^{q}$ a aucune représentation. Donc notre conjoncture que $\left(h_{n}\right)$ existe est fausse et $\left(g_{n}\right)$ est unique.

\section{Conclusion}
 La Généralisation du Théorème de Zeckendorf  ne doit pas rester plus qu’une curiosité. Dans le futur
recherche, nous prévoyons d’étudier les applications de nos résultats à d’autres domaines de
mathématiques telles que les codes correcteurs d’erreurs.

\medskip



\newpage
\begin{thebibliography}{99}
\bibitem{1} J.L Broun,\emph{Unique Representation of integers as sums of Distinct Lucas Numbers},
The Fibonacci Quarterly, Vol 7, No 3, (1969) pp 243-252.
\bibitem{2} E. Zeckendorf,\emph{ Représentation des nombres naturels par une somme de nombres de Fibonacci ou de nombres de Lucas} ,Bull. Soc. R. Sci. Liège Vol. 41 (1972) pp. 179 – 182.
\bibitem{3} D. E. Daykin,\emph { Representation of Natural Numbers as Sums of Generalized Fibonacci Numbers,} J.London Math .soc. , 35 (I960), pp. 143-160.
 
\bibitem{4} P. Filiponi,\emph{The Representation of Certain Integers as a Sum of Distinct
Fibonacci Numbers}, Tech Rep. 2B0985. Fondazione Ugo Bordoni, Rome (1985).
\bibitem{5} Rachid Chergui,\emph{Zeckendorf Arithmetic For Lucas Numbers},Palestine Journal of Mathematics, Vol. 9(1)(2020) , 337–342.
\bibitem{6} R.L. Graham,, D.E. Knuth,, and O, Patashnik,\emph{ Concrete Mathematics}, Addison-Wesley, Reading, MA 1991, pp. 295-297.






\end{thebibliography}
\end{document}